\documentstyle{article}
\begin{document}
\newcommand{\HDS}{\vrule width0pt height2.3ex depth1.05ex\displaystyle}
\def\f#1#2{{{\HDS #1}\over{\HDS #2}}}
\def\lpravilo#1{ \makebox[-.5em][r]{\mbox{\it #1}} {\mbox{\hspace{0.5em}}}}

\title{\bf Identity of Proofs \\Based on \\Normalization and Generality}

\author{{\sc Kosta Do\v sen}
\\[.05cm]
\\Mathematical Institute, SANU \\
Knez Mihailova 35, p.f. 367 \\
11001 Belgrade, Serbia \\
email: kosta@mi.sanu.ac.yu}
\date{}
\maketitle

\begin{abstract}
\noindent Some thirty years ago, two proposals were made concerning criteria
for identity of proofs. Prawitz proposed to analyze identity of proofs in
terms of the equivalence relation based on reduction to normal form in
natural deduction. Lambek worked on a normalization proposal analogous to
Prawitz's, based on reduction to cut-free form in sequent systems, but he
also suggested understanding identity of proofs in terms of an equivalence
relation based on generality, two derivations having the same generality if
after generalizing maximally the rules involved in them they yield the same
premises and conclusions up to renaming of variables. These two proposals
proved to be extensionally equivalent only for limited fragments of logic.

The normalization proposal stands behind very successful applications of the
typed lambda calculus and of category theory in the proof theory of
intuitionistic logic. In classical logic, however, it did not fare well.

The generality proposal was rather neglected in logic, though related
matters were much studied in pure category theory in connection with
coherence problems, and there are also links to low-dimensional topology and
linear algebra. This proposal seems more promising than the other one for
the general proof theory of classical logic.

\vspace{0.3cm}

\noindent{\it Mathematics Subject Classification} ({\it 2000}): 03F03, 03F07,
03A05, 03-03 \\[.2cm]
{\it Keywords}: proof, criteria of identity, cut elimination and normal
form, generality, categorial coherence
\end{abstract}

\section{General proof theory}

In mathematics, as in other established bodies of human knowledge, when a
question is not in the main stream of investigations, it runs the danger
of being given short shrift by the establishment. Often, when
they are of a certain general kind, such
questions are dismissed as being ``philosophical'', and often this need
not mean that the question is inherently such, that it
cannot be approached in a mathematical spirit, but it means only that
the person dismissing the question does not want to deal with it,
presumably because he does not know how.

Before the advent of recursion theory, many mathematicians would
presumably have dismissed the question ``What is a computable function?''
as a philosophical question, unworthy of attention, which only an
outsider would deal with.
(I suppose that this question could even have been understood as a
psychological, empirical, question.)
It required something like the enthusiasm of a
young discipline on the rise, which logic was in the first half of the
twentieth century, for such questions to be embraced as legitimate, and
seriously treated by mathematical means---with excellent results. In the
meantime, logic has aged, with losses in enthusiasm and gains in
conformism.

An outsider might suppose that the question ``What is a proof?'' should be
important for a field called {\it proof theory}, and then he
would be surprised to find that this 
and related questions, one of which will occupy us here, are
exactly of the kind to be dismissed as ``philosophical'' by the
establishment. These questions are of the {\it conceptual}
kind, and not of the {\it technical} kind, but a great deal of
mathematics---a great deal of the best mathematics---is of the
conceptual kind.

By the beginning of the 1970s, Prawitz, inspired by ideas of Kreisel,
inaugurated the field of {\it general proof theory} with this and related
questions (see \cite{P71}, Section I.2; see also the next section of this
survey).
The field exists, I believe, but its borders are not well
marked, and the name Prawitz gave to it is not widespread.
(Some authors speak nowadays of {\it structural} proof theory, but
I am not sure this kind of proof theory, which seems to amount to
{\it basic} proof theory, is Prawitz's general proof theory.)
We need,
however, a name for that part of proof theory whose foundations are in
Gentzen's thesis \cite{G35}, and which differs from the proof theory 
whose subject grew out of consistency proofs for formalized
fragments of mathematics inaugurated also by Gentzen.
This other kind of proof theory, whose foundations are
in Hilbert's program and G\" odel's results, and 
which was dominant throughout
the twentieth century, was called {\it reductive} by Prawitz.

General proof theory, although it is seemingly more philosophical,
is closer to applications than reductive proof theory. It is in the period
when logic became increasingly intertwined with computer science, in
the last decades of the twentieth century, that general proof theory started
gaining some ground. Two matters mark this relative, still limited,
success: one is the connection with the
typed lambda calculus, based on the Curry-Howard correspondence, and
the other
is the connection with category theory. A third area of proof theory,
distinct from general proof theory and reductive proof theory, emerged
at roughly the same period in association with complexity theory
(see \cite{P98}).
Hilbert's name is usually tied to reductive proof theory,
but from the very beginning of the century
he was also interested in some questions of
general proof theory that touch upon complexity theory
(see \cite{T03}).

\section{The Normalization Conjecture}

To keep up with the tradition, throughout the present survey we speak of
``proof'', though we could as well 
replace this term by the term  ``deduction'', which might
be more precise, since we have in mind deductive proofs from assumptions
(including the empty collection of assumptions).

The question ``What is a proof?'', as the basic question of the field of
general proof theory, was introduced by Prawitz in \cite{P71}
(Section I) with
additional specifications. Let us quote the relevant passage in
extenso:

\begin{quotation}
\noindent Obvious topics in general proof theory are:

\vspace{.1cm}

\noindent 2.1. The basic question of defining the notion of proof, including
the question of the distinction between different kinds of proofs
such as constructive proofs and classical proofs.

\vspace{.1cm}

\noindent 2.2. Investigation of the structure of (different kinds of) proofs,
including e.g. questions concerning the existence of certain normal
forms.

\vspace{.1cm}

\noindent 2.3. The representation of proofs by formal derivations.
In the same way
as one asks when two formulas define the same set or two sentences
express the same proposition, one asks when two derivations represent
the same proof; in other words, one asks for identity criteria for proofs or
for a ``synonymity'' (or equivalence) relation between derivations.

\vspace{.1cm}

\noindent 2.4. Applications of insights about the structure of proofs
to other
logical questions that are not formulated in terms of the notion of proof.
(p. 237)
\end{quotation}

We are here especially interested in Prawitz's topic 2.3, namely, in the
question of identity criteria for proofs. Prawitz's merit is that he
did not only formulate this question very clearly, but also proposed a
precise mathematical answer to it.

Prawitz considered derivations in natural deduction systems and the
equivalence relation between derivations that is the reflexive, transitive
and symmetric closure of the immediate reducibility relation between
derivations. Of course, only derivations with the same assumptions
and the same conclusion may be equivalent. Prawitz's immediate
reducibility relation is the one involved in reducing a derivation to
normal form.
(As it is well known, the idea of this reduction stems
from Gentzen's thesis \cite{G35}.)
A derivation reduces immediately to another derivation
(see \cite{P71}, Section II.3.3) when the latter is obtained from the former
either by removing a {\it maximum formula} (i.e.\ a formula with a connective
$\alpha$ that is the conclusion of an introduction of $\alpha$ and the
major premise of an elimination of $\alpha$), or by performing one
of the {\it permutative reductions} tied to the eliminations of
disjunction and of the existential quantifier, which enables us to
remove what Prawitz calls {\it maximum segments}. There are some
further reductions, which Prawitz called {\it immediate simplifications};
they consist in removing eliminations of disjunction where no
hypothesis is discharged, and there are similar immediate
simplifications involving the existential quantifier, and
``redundant'' applications of the classical absurdity rule. Prawitz
also envisaged reductions he called {\it immediate expansions}, which
lead to the expanded normal form where all the {\it minimum formulae}
are atomic (minimum formulae are those that are conclusions
of eliminations and premises of introductions).

Prawitz formulates in \cite{P71} (Section II.3.5.6) the following
conjecture, for which he gives credit (in Section II.5.2) to
Martin-L\" of, and acknowledges influence by ideas of Tait:

\begin{quotation}
\noindent {\bf Conjecture.} Two derivations represent the same proof
if and only if they are equivalent.
\end{quotation}

\noindent We will call this conjecture the {\it Normalization
Conjecture}.

Prawitz considered the right-to-left direction of the Normalization Conjecture
as relatively unproblematic. He found it more difficult 
to unearth facts that would support the other direction of the conjecture.
At the same time (in the same book where Prawitz's paper
\cite{P71} appeared), Kreisel
considered how one could justify by mathematical means a conjecture
such as the Normalization Conjecture, and in particular the left-to-right
direction of the conjecture (see \cite{K71}, pp.\ 114-117, 165).

The Normalization Conjecture is an assertion of the same kind as Church's
Thesis: we should not expect a formal proof of it. (Kreisel speaks in
\cite{K71}, p. 112, of the ``informal'' character of the conjecture.)
The Normalization Conjecture attempts to give a formal reconstruction
of an intuitive notion. (Like Church's Thesis, the Normalization Conjecture
might be taken as a kind of definition. It is, however, better to
distinguish this particular kind of definition by a special name.
The Normalization Conjecture, as well as the Generality Conjecture,
which we will consider in the next section, might be taken as a case of
{\it analysis} in the sense of \cite{D89}.)

Kreisel suggests (in \cite{K71}, p. 116) that the right-to-left direction of
the Normalization Conjecture may be seen as analogous to a {\it soundness}
assertion, while the left-to-right direction is analogous to a
{\it completeness} assertion. Derivations are the syntactical matter,
while intuitive proofs are their semantics. One could, however, take
an exactly opposite view: proofs are the given syntactical matter,
which we model by derivations.

Martin-L\" of examined the
Normalization Conjecture in \cite{Mar75a} (Section 2.4, p. 104),
giving credit to Prawitz for its formulation.
He maintained that there is little hope of establishing this conjecture
unless
``represent the same proof'' on the left-hand side is replaced by
``are definitionally equal'', or identity of proofs on the same side
is replaced by {\it provable} identity; in both cases, the matter is
treated further within the realm of his intuitionistic theory of
types of \cite{Mar75}.

The Normalization Conjecture was formulated by Prawitz at the time when
the Curry-Howard correspondence between derivations in natural deduction and
typed lambda terms started being recognized more and more (though the
epithet ``Curry-Howard'' was not yet canonized). Prawitz's equivalence
relation between derivations corresponds to beta-eta equality
between typed lambda terms, if immediate expansions are taken into
account, and to beta equality otherwise.

When speaking of the Curry-Howard correspondence, one should, however, bear in
mind that it covers really well only the conjunction-implication fragment
of intuitionistic propositional logic. In the presence of disjunction, the
typed lambda calculus is extended with a ternary variable-binding
operation that serves to code disjunction elimination. Such lambda
calculuses exist, but they are not particularly appealing. Negation, i.e.\
the absurd constant proposition, brings its own complications.
(With variables of the type of the absurd constant proposition
one derives every equation.)

Nevertheless, the fact that we have an alternative formal representation of
proofs with typed lambda terms, and that equality between these terms
agrees so well with Prawitz's equivalence of derivations, is quite
remarkable. This fact lends support to the Normalization
Conjecture, as the fact that alternative formal representations of
computable functions cover the same functions of arithmetic lends
support to Church's Thesis.

Besides derivations in natural deduction
and typed lambda terms, where according to the Curry-Howard
correspondence the latter can be conceived just as codes for the former,
there are other, more remote, formal representations of proofs.
There are first
Gentzen's sequent systems, which are related to natural
deduction, but are nevertheless different, and there are also
representations of proofs as arrows in categories. Equality of arrows
with the same domain and codomain, i.e.\ commuting diagrams of arrows,
which is what category theory is about, should now correspond to
identity of proofs via a conjecture analogous to the Normalization
Conjecture.

The fact proved by Lambek (see \cite{L74} and \cite{LS86}, Part I;
see also \cite{D01}) that the category of typed lambda calculuses
with functional types and finite product types, based on
beta-eta equality, is equivalent to the category
of cartesian closed categories, and that hence equality of typed
lambda terms amounts to equality between arrows in cartesian
closed categories, lends additional support to the Normalization
Conjecture. Equality of arrows in bicartesian closed categories
(i.e.\ cartesian closed categories with finite coproducts)
corresponds to equivalence of derivations in Prawitz's sense in 
full intuitionistic propositional logic.

In category theory, the Normalization Conjecture is tied to Lawvere's
Thesis that all the logical constants of intuitionistic logic are
characterized by adjoint situations. Prawitz's equivalence of
derivations, in its beta-eta version, corresponds to equality of
arrows in various adjunctions tied to logical constants (see
\cite{Law69}, \cite{D99}, Section 0.3.3, \cite{D99a} and \cite{D01}).
Adjunction is the unifying concept for the reductions envisaged
by Prawitz.

Even the fact that equality between lambda terms, as well as equality
of arrows in cartesian closed categories, were first conceived for
reasons independent of proofs is remarkable. This tells us that we are
in the presence of a solid mathematical structure, which may be
illuminated from many sides, and is not a figment.

Prawitz formulated the Normalization Conjecture having in mind
mainly intuitionistic logic. It was meant to cover classical logic too.
One may, however, formulate an analogous conjecture for all logics
with a normalization procedure analogous to Prawitz's procedure
for intuitionistic logic. Such are in particular substructural logics.
Such are also modal extensions of various logics. It is not evident
that Lawvere's Thesis applies to all the logical constants of all these
logics; the adjunctions in question may be hidden, or perhaps there
is no adjunction, but only something related to it (for example,
something like an adjunction but lacking the unit or the
counit natural transformation).

It may be asked whether the Normalization Conjecture, or another
answer to the
question ``When are two representations of proofs equal?'', brings us
closer to answering the basic question ``What is a proof?''. If our
answers to the first question about identity criteria for proofs cut
across various representations of proofs, so that when one representation is
translated into another they agree about the notion of identity of proof,
it seems that we could answer the basic question ``What is a proof?'' by
taking one of our representations and declaring that the equivalence classes
of this representation determine the notion we are seeking. A proof is
the equivalence class of one of its representations. This does not
settle, however, which representation we should choose. It presumably
cannot be any representation, but one that neither abolishes important
distinctions in the structure of the proof, nor introduces irrelevant
features. Anyway, it seems the Normalization Conjecture brings us closer to
answering the basic question.

The answer to the question ``What is a computable function?'', which is
given by Church's Thesis, did not involve such a factoring through an
equivalence relation, but the related question ``What is an algorithm?''
may require exactly that. This last question was raised in recursion
theory fairly recently by Moschovakis in \cite{M01}, where the
fundamental interest of identity criteria for algorithms is stressed
(see Section 8 of that paper). It is rather natural to find that the
questions ``What is a proof?'' and ``What is an algorithm?'' are
analogous. As a matter of fact, many logicians conceive proofs just
as a kind of algorithm for generating theorems. What should
distinguish, however, this particular kind of algorithm is the
role logical constants play in it.

While Moschovakis approaches the question of identity criteria for
algorithms in a mathematical spirit, proof-theorists nowadays
most often ignore the question of identity criteria for proofs, or
they dismiss it as ``philosophical'', sometimes expressing doubt that a
simple answer will ever be given to this question. The answer
Prawitz gave to this question with the Normalization Conjecture
is certainly not philosophical if
``philosophical'' means ``vague'' or ``imprecise''. Prawitz's answer,
which is rather simple, is maybe not conclusive, but it is a mathematical
answer. All the more so if we take into account its connection with the
lambda calculus and fundamental structures of category theory.

\section{The Generality Conjecture}

At the same time when Prawitz formulated the Normalization Conjecture,
in a series of papers (\cite{L68}, \cite{L69}, \cite{L72} and \cite{L74})
Lambek was engaged
in a project where arrows in various sorts of
categories were construed as representing proofs. The domain of an arrow
corresponds to the collection of assumptions joined by a conjunction
connective, and the codomain corresponds to the conclusion. With this series
of papers Lambek inaugurated the field of categorial proof theory.

The categories Lambek considered in \cite{L68} and \cite{L69} are first
those that correspond to his substructural syntactic calculus of
categorial grammar (these are monoidal categories where the functors
$A\otimes \ldots$ and $\ldots \otimes A$ have right adjoints). Next,
he considered monads, which besides being fundamental for category theory,
cover proofs in modal logics of the S4 kind. In \cite{L72} and \cite{L74}
Lambek dealt with cartesian closed categories, which is Prawitz's ground,
since these categories cover proofs in the conjunction-implication fragment
of intuitionistic logic.  He also envisaged bicartesian closed categories,
which cover the whole of intuitionistic propositional logic.

Lambek's insight is that equations between arrows in categories, i.e.\
commuting diagrams of arrows, guarantee cut elimination, i.e.\
composition elimination, in an appropriate language for naming arrows.
(In \cite{D99} it is established that for some basic notions
of category theory, and in particular for the notion of adjunction,
the equations assumed are necessary and sufficient for cut
elimination.)
Since cut elimination is closely related to Prawitz's normalization of
derivations, the equivalence relation envisaged by Lambek should be related
to Prawitz's.

Actually, Lambek dealt with cut elimination in categories related to
his syntactic calculus and in monads. For cartesian closed categories he
proved in \cite{L74} another sort of result, which he called
{\it functional completeness}. Lambek's functional completeness is
analogous to the functional completeness (often called {\it combinatory
completeness}) of Sch\" onfinkel's and Curry's systems of combinators,
which enables us to define lambda abstraction in these systems. With
the help of functional completeness, via the typed lambda calculus,
Lambek established that his equivalence relation between derivations
in conjunctive-implicative intuitionistic logic and Prawitz's
beta-eta equivalence amount to the same thing, as we mentioned in
the preceding section.

Lambek's work is interesting for us here not only because 
he worked with an equivalence relation between derivations amounting to
Prawitz's, but also because he envisaged another kind of equivalence
relation. Lambek's idea is best conveyed by considering the following example.
In \cite{L72} (p. 65) he says that the first projection arrow
$\pi^1_{p,p}: p\wedge p \rightarrow p$ and the second projection arrow
$\pi^2_{p,p}: p\wedge p \rightarrow p$, which correspond to
two derivations of conjunction elimination,
have different {\it generality},
because they generalize to $\pi^1_{p,q}: p\wedge q \rightarrow p$ and
$\pi^2_{p,q}: p\wedge q \rightarrow q$ respectively, and these two arrows
do not have the same codomain; on the other hand,
$\pi^1_{p,q}: p\wedge q \rightarrow p$ and
$\pi^2_{q,p}: q\wedge p \rightarrow p$ do not have the same domain.
The idea of generality may be explained roughly as follows.
We consider generalizations of derivations that diversify variables
without changing the rules of inference.
Two derivations have the same generality when every
generalization of one of them
leads to a generalization of the other, so that in the
two generalizations we have the same assumptions and conclusion
(see \cite{L68}, p. 257). In the example above this is not the case.

Generality induces an equivalence relation between derivations.
Two derivations are equivalent if and only if they have the same
generality. Lambek does not formulate so clearly as Prawitz a conjecture
concerning identity criteria for proofs, but he suggests 
that two derivations represent the same proof if and only if they are
equivalent in the new sense. We will call this conjecture the
{\it Generality Conjecture}.

The left-to-right direction of the Generality Conjecture seems
pretty intuitive. The other direction is less clear, and if it
happens to be true, that would be something like a discovery.
With the Normalization Conjecture, the situation was different:
the right-to-left direction seemed unproblematic, and the
left-to-right direction was less clear.

Lambek's own attempts at making the notion of generality precise
(see \cite{L68}, p. 316, where the term ``scope'' is used instead
of ``generality'', and \cite{L69}, pp.\ 89, 100) need not detain us here.
In \cite{L72} (p. 65) he finds that these attempts were faulty.

The simplest way to understand generality is to use graphs whose
vertices are occurrences of propositional letters in the assumptions and the
conclusion of a derivation. We connect by an edge occurrences of 
letters that must remain occurrences of the same letter after generalizing,
and do not connect those that may become occurrences of different letters.
So for the first and second projection above we would have the two graphs

\begin{center}
\begin{picture}(180,50)
\put(40,40){\makebox(0,0){$\wedge$}}
\put(155,40){\makebox(0,0){$\wedge$}}
\put(25,40){\makebox(0,0){$p$}}
\put(0,25){\makebox(0,0){$\pi^1_{p,p}$}}
\put(55,40){\makebox(0,0){$p$}}
\put(140,40){\makebox(0,0){$p$}}
\put(115,25){\makebox(0,0){$\pi^2_{p,p}$}}
\put(170,40){\makebox(0,0){$p$}}
\put(40,10){\makebox(0,0){$p$}}
\put(155,10){\makebox(0,0){$p$}}
\put(27,30){\line(1,-1){12}}
\put(168,30){\line(-1,-1){12}}
\end{picture}
\end{center}

\noindent When the propositional letter $p$ is replaced by an
arbitrary formula
$A$ we have an edge for each occurrence of propositional letter in $A$.

The generality of a derivation is such a graph. According to the
Generality Conjecture, the first and second projection derivations
from $p \wedge p$ to $p$ represent different proofs because their
generalities differ.

One defines an associative composition of such graphs, and there is also
an obvious identity graph with straight parallel edges, so that
graphs make a category, which we call the {\it graphical category}.
If on the other hand it is taken for granted that proofs
also make a category, which we will call the {\it syntactical category},
with composition of arrows being composition of
proofs, and identity arrows being identity proofs
(an identity proof composed with any
other proof, either on the side of the assumptions or on the side
of the conclusion, is equal to this other proof),
then the Generality Conjecture may be rephrased as the assertion that
there is a faithful functor from the syntactical category to
the graphical category.

If our syntactical category is precisely determined in advance,
by relying on the Normalization Conjecture, or perhaps in another
way, then the Generality Conjecture is something we
can prove or disprove formally. If on the other hand this
syntactical category is left undetermined, then the Generality Conjecture
is in the same position as the Normalization Conjecture and Church's
Thesis, and cannot be proved or disproved formally. We can only
try to see how well it accords with our other intuitions, and in
particular with the intuition underlying the Normalization Conjecture.

For example, it may be taken that 
the syntactical category in conjunctive logic (without implication), both
intuitionistic and classical (these two conjunctive logics do not differ),
is a cartesian category $\cal K$ freely
generated by a set of propositional letters (a category is cartesian when
it has all finite products). Conjunctive logic here includes the
true constant proposition. That this syntactical category may indeed
be taken as the category of proofs in conjunctive logic is
in accordance with the Normalization Conjecture in its beta-eta version,
but it is not excluded that we have chosen this category for other
reasons. 
The graphical category $\cal G$ in this
particular case may be taken to be the category opposite to the category
of functions on finite ordinals. The finite ordinals correspond to the
number of occurrences of propositional letters in a formula, and functions go
from the conclusions to the assumptions; that is why we have the opposite
of the category of functions. The Generality Conjecture is the assertion
that there is a faithful functor from $\cal K$ to $\cal G$; that is, a
functor $G$ from $\cal K$ to $\cal G$ such that for $f$ and $g$ with the 
same domain and the same codomain we have

\begin{center}
$(*)$ \quad \quad $f=g$ {\it in} $\cal K \;\;$ {\it if and only if} 
$\;\;\; G(f)=G(g)$ {\it in} $\cal G$.
\end{center}

\noindent From left to right, the equivalence $(*)$ follows from the
functoriality of $G$, and from right to left, it expresses the
faithfulness of $G$. In this particular case, the Generality Conjecture
understood as $(*)$ can be proved formally
(see \cite{P02}, \cite{DP01} and references
therein concerning earlier statements and proofs of equivalent
results).

Note that graphs are not representations of proofs, as derivations in natural
deduction were. A graph may abolish many important distinctions; in
particular, the images of different connectives, like conjunction and
disjunction, may be isomorphic in the graph.
The Generality Conjecture answers the question of identity
criteria for proofs, and brings us closer to answering the basic
question ``What is a proof?'', by involving
representations of proofs more plausible than graphs, like those
we find in natural deduction, or in sequent systems, or in
the syntactical categories. A proof would
be the equivalence class of such a plausible representation with
respect to the equivalence relation defined in terms of generality,
i.e.\ in terms of graphs.

Understood in the sense of $(*)$, the Generality Conjecture is analogous to
what in category theory is called a {\it coherence theorem}. Coherence in
category theory was understood in various ways (we cannot survey here the
literature on this question; for earlier works see \cite{ML72}), but
the paradigmatic coherence results of \cite{ML63} and \cite{KML71} can
be understood as faithfulness results like $(*)$, and this is how we will
understand coherence. (The coherence questions one encounters in
the theory of weak {\it n}-categories, which we will mention
at the end of the
next section, are of another, specific, kind.)

The coherence result of \cite{KML71} proves the Generality Conjecture
for the multiplicative conjunction-implication fragment of 
intuitionistic linear logic (modulo a condition concerning the
multiplicative true constant proposition, i.e.\ the unit with
respect to multiplicative conjunction), and, inspired by Lambek, it
does so via a cut-elimination proof. The syntactical category in this case
is a free symmetric monoidal closed category, and the graphical category
is of a kind studied in \cite{EK66}. The graphs of this graphical
category are closely related to what in knot theory is called
{\it tangles}. In tangles, as in braids, we distinguish between two
kinds of crossings, but here we need just one kind, in which it is not
distinguished which of the two crossed edges is above the other.
(For categories of tangles see \cite{Y88}, \cite{T89} and \cite{K95},
Chapter 12.) Tangles with this single kind of crossing are like
graphs one encounters in Brauer algebras (see \cite{B37} and \cite{W88}).
Here is an example of such a tangle:

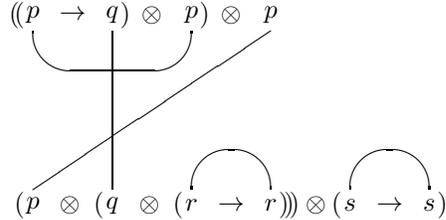
\begin{center}
\begin{picture}(130,95)

\put(0,15){\line(3,2){90}}
\put(30,15){\line(0,1){60}}

\put(75,15){\oval(30,30)[t]}
\put(135,15){\oval(30,30)[t]}
\put(30,75){\oval(60,30)[b]}

\put(-5,9){\makebox(0,0){$($}}
\put(0,9){\makebox(0,0){$p$}}
\put(14,9){\makebox(0,0){$\otimes$}}
\put(25,9){\makebox(0,0){$($}}
\put(30,9){\makebox(0,0){$q$}}
\put(44,9){\makebox(0,0){$\otimes$}}
\put(55,9){\makebox(0,0){$($}}
\put(60,9){\makebox(0,0){$r$}}
\put(75,9){\makebox(0,0){$\rightarrow$}}
\put(90,9){\makebox(0,0){$r$}}
\put(95,9){\makebox(0,0){$)$}}
\put(97,9){\makebox(0,0){$)$}}
\put(99,9){\makebox(0,0){$)$}}
\put(107,9){\makebox(0,0){$\otimes$}}
\put(115,9){\makebox(0,0){$($}}
\put(120,9){\makebox(0,0){$s$}}
\put(135,9){\makebox(0,0){$\rightarrow$}}
\put(150,9){\makebox(0,0){$s$}}
\put(155,9){\makebox(0,0){$)$}}

\put(-7,81){\makebox(0,0){$($}}
\put(-5,81){\makebox(0,0){$($}}
\put(0,81){\makebox(0,0){$p$}}
\put(15,81){\makebox(0,0){$\rightarrow$}}
\put(30,81){\makebox(0,0){$q$}}
\put(35,81){\makebox(0,0){$)$}}
\put(45,81){\makebox(0,0){$\otimes$}}
\put(60,81){\makebox(0,0){$p$}}
\put(65,81){\makebox(0,0){$)$}}
\put(75,81){\makebox(0,0){$\otimes$}}
\put(90,81){\makebox(0,0){$p$}}

\end{picture}
\end{center}

Tangles without crossings at all serve in \cite{D99} (Section 4.10;
see also \cite{D03})
to obtain a coherence result for the general notion of adjunction,
which according to Lawvere's Thesis underlies all logical constants
of intuitionistic logic, as we mentioned in the preceding section.
In terms of combinatorial low-dimensional topology, the mathematical
content of the general notion of adjunction is caught by the
Reidemeister moves of planar ambient isotopy. An analogous
coherence result for self-adjunctions, where a single endofunctor
is adjoint to itself, is proved in \cite{DP02a}. Through this 
latter result we reach the theory of Temperley-Lieb algebras,
which play a prominent role in knot theory and low-dimensional
topology, due to Jones' representation of Artin's braid groups in
these algebras (see \cite{KL}, \cite{L97}, \cite{PS}
and references therein).

In \cite{DP02a} 
one finds also coherence results for self-adjunctions
where the graphical category is the category of matrices,
i.e.\ the skeleton of the category of finite-dimensional
vector spaces over a fixed field with linear transformations as arrows.
Tangles without crossings may be faithfully represented in matrices by
a representation derived from
the orthogonal group case of Brauer's representation of Brauer
algebras (see also \cite{W88}, Section 3, and \cite{J94}, Section 3).
This representation is based on the fact that the Kronecker product of
matrices gives rise to a self-adjoint functor in the category of matrices,
and this self-adjointness is related to the fact that in this category,
as well as in the category of binary relations between finite ordinals,
finite products and coproducts are isomorphic.

In \cite{P02} there are interesting coherence results,
which extend \cite{ML63}, for the multiplicative conjunction
fragments of substructural logics. Further coherence results for
intuitionistic and classical logic will be considered in the next
section.

In the light of the connection with combinatorial aspects of low-dimensional
topology, which we mentioned above, the following quotation from
\cite{K71} sounds prophetical:

\begin{quotation}
\noindent {\it Remark.} At the risk of trying to explain {\it obscurum
per obscurius}, I should point out a striking analogy between the
problem of finding $(\alpha)$ a conversion relation corresponding to
identity of proofs and $(\beta)$ an ``equivalence'' relation such as
that of {\it combinatorial equivalence} in topology corresponding to
a basic invariant in our geometric concepts. (p. 117)
\end{quotation}

\noindent (The sequel of Kreisel's text suggests he had in mind
equivalence relations that cut across various representations
of proofs, i.e.\ are insensitive to the peculiarities of
particular formalizations.)

From a logical point of view, the Generality Conjecture may be understood as
a completeness theorem. The left-to-right direction of the equivalence
$(*)$ is
{\it soundness}, and the right-to-left direction is {\it completeness}
proper. (This terminology is opposite to that suggested by Kreisel
for the Normalization Conjecture, which we mentioned in the preceding
section.)

Like the Normalization Conjecture, the Generality Conjecture gives a precise
mathematical answer to the question about identity criteria for proofs, and,
like the former conjecture, it has ties with important mathematical structures.

\section{The two conjectures compared}

The Normalization Conjecture and the Generality Conjecture agree
only for limited fragments of logic. As we said in the preceding section,
they agree for purely
conjunctive logic (without implication),
with or without the true constant proposition
$\top$. Conjunctive logic is the same for intuitionistic and classical
logic. Here the Normalization Conjecture is taken in its beta-eta
version. By duality, the two conjectures agree for purely disjunctive
logic, with or without the absurd constant proposition $\bot$. If we
have both conjunction and disjunction, but do not yet have distribution
of conjunction over disjunction, and have neither $\top$ nor $\bot$, then the
two conjectures still agree (see \cite{DP02}). And here it seems we
have reached the limits of agreement as far as intuitionistic and
classical logic are concerned, provided the graphs
involved in the Generality Conjecture correspond to relations between
the domain and the codomain.
With more sophisticated notions of graphs, matters may stand
differently, and the area of agreement for the two conjectures 
may perhaps be wider, but it can be even narrower, as we will see below.

Coherence fails for the syntactical
bicartesian categories, in which besides nonempty
finite products and coproducts we have also an empty product and
and an empty coproduct, i.e.\ a terminal object $\top$ and an
initial object $\bot$
(an object in a category is terminal when there is exactly one arrow
from any object to it, and it is initial when there is exactly one arrow
from it to any object). In order to regain coherence we must pass to special
bicartesian categories where the first projection from the product of
$\bot$ with itself $\pi^1_{\bot , \bot}: \bot \times \bot \rightarrow \bot$
is equal to the second projection $\pi^2_{\bot , \bot}: \bot
\times \bot \rightarrow \bot$, and analogously for the first and second
injection from $\top$ to the coproduct of $\top$ with itself, namely
$\top + \top$ (see \cite{DP01a}). Here $\times$ corresponds to conjunction
and $+$ to disjunction. The assumption concerning these two projections
holds in the category {\bf Set}, but the assumption concerning the two
injections does not. Both assumptions hold in some other important
bicartesian categories: the category of pointed sets, and its subcategories
of commutative monoids with monoid homomorphisms and of vector spaces
over a fixed field with linear transformations. They hold also in 
the category of matrices mentioned
in the preceding section. Anyway, in the presence of $\top$ and $\bot$
the two conjectures do not agree, since normalization does not deliver
the assumption concerning the two injections.

The graphical category with respect to which
we have coherence for syntactical categories
with nonempty finite products and
coproducts is the category of binary relations between finite ordinals
(see \cite{DP02}). In that category, $\times$ and $+$, which correspond
to conjunction and disjunction respectively, are isomorphic.

The intuitive idea of generality for conjunctive-disjunctive logic
is not caught by arbitrary binary relations, but by binary relations that
are {\it difunctional} in the sense of \cite{R48} (Section 7). A binary
relation $R$ is difunctional when $RR^{-1}R \subseteq R$. The
composition of two difunctional relations is not necessarily
difunctional, so that we do not have at our disposal the category
of difunctional relations, with respect to which we could prove
coherence, and, anyway, the image under the
faithful functor of an arrow
from our syntactical category with nonempty finite products and
coproducts is not necessarily a difunctional relation.

Even if difunctionality were satisfied, it could still be questioned
whether our intuitive idea of generality is caught by binary relations
in the case of conjunctive-disjunctive logic. The problem is that if
$w_p:p \rightarrow p \times p$ is a component of the diagonal
natural transformation, and $\iota^1_{q,p}:q \rightarrow q+p$ is a first
injection, then in categories with finite products and coproducts we have
\[
({\mbox{\bf 1}}_q + w_p) \circ \iota^1_{q,p} = \iota^1_{q,p \times p}
\]
where the left-hand side cannot be further generalized, but the right-hand
side can be generalized to $\iota^1_{q, p \times r}$.
The intuitive idea of generality seems to require that in
$w_p:p \rightarrow p \times p$ we should not have only a relation
between the domain and the codomain, as on the left-hand side below,
but a graph as on the right-hand side:

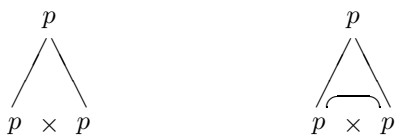
\begin{center}
\begin{picture}(193,50)
\put(40,5){\makebox(0,0){$\times$}}
\put(155,5){\makebox(0,0){$\times$}}

\put(27,5){\makebox(0,0){$p$}}
\put(53,5){\makebox(0,0){$p$}}
\put(142,5){\makebox(0,0){$p$}}
\put(168,5){\makebox(0,0){$p$}}

\put(40,45){\makebox(0,0){$p$}}
\put(155,45){\makebox(0,0){$p$}}

\put(26,12){\line(1,2){13}}
\put(54,12){\line(-1,2){13}}
\put(141,12){\line(1,2){13}}
\put(169,12){\line(-1,2){13}}
\put(155,12){\oval(20,8)[t]}

\end{picture}
\end{center}

\noindent (see \cite{DP021}, and also \cite{DP022}). With such graphs
we can still get coherence for
conjunctive logic, and for disjunctive logic, taken separately, but for
conjunctive-disjunctive logic the left-to-right direction, i.e.\
the soundness part, of coherence would fail. So for conjunctive-disjunctive
logic the idea of generality with which we have coherence is not quite
the intuitive idea suggested by Lambek, but only something close to it,
which involves the categorial notion of natural transformation.

The fragments of logic mentioned above where the Normalization
Conjecture and the Generality Conjecture agree all possess a
property called {\it maximality}. Let us say a few words about this
very important property.

For the whole field of general proof theory to make sense, and in
particular for considering the question of identity criteria for
proofs, we should not have that any two derivations with the same
assumptions and conclusion are equivalent, i.e.\, it should not be the case
that there is never more than one proof with given assumptions and
a given conclusion. Otherwise, our field would be trivial.

This marks the watershed between proof theory and the rest of logic, where
one is not concerned with proofs, but at most with consequence
{\it relations}. With relations, we either have a pair made of
a collection of assumptions and a conclusion, or we do not have it.
In proof theory, such pairs are indexed with various
proofs, and there may be several proofs for a single pair.

Now, categories with finite nonempty products, cartesian categories
and categories with finite nonempty products and coproducts
have the following property. Take, for example, cartesian categories,
and take any equation in the language of free cartesian categories that
does not hold in free cartesian categories. If a cartesian category
$\cal K$ satisfies this equation, then $\cal K$ is a {\it preorder};
namely, all arrows with the same domain and codomain are equal.
We have an exactly analogous property with the other sorts of
categories we mentioned. This property is a kind of Post completeness.
Any extension of the equations postulated leads to collapse.

Translated into logical language, this means that Prawitz's
equivalence relation for derivations in conjunctive logic,
disjunctive logic and conjunctive-disjunctive logic without
distribution and without $\top$ and $\bot$, which in all these
cases agrees with our equivalence relation defined via
generality in the style of Lambek, is {\it maximal}.
Any strengthening, any addition, would yield
that any two derivations with the same assumptions and the same
conclusion are equivalent.

If the right-to-left direction of the Normalization Conjecture holds,
with maximality we can efficiently justify the left-to-right direction,
which Prawitz found problematic
in \cite{P71}, and about which Kreisel was thinking in \cite{K71}.
In the footnote on p. 165 of that paper Kreisel mentions 
that Barendregt suggested this justification via maximality. 
Suppose the right-to-left direction of the Normalization Conjecture
holds, suppose that for some assumptions and conclusion there is more
than one proof, and suppose the equivalence relation is maximal. 
Then if two derivations represent the same proof, they are equivalent.
Because if they were not equivalent, we would never have more
than one proof with given assumptions and a given conclusion.
Nothing can be missing from our equivalence relation, because whatever
is missing, by maximality, leads to collapse on the side of the
equivalence relation, and, by the right-to-left direction of
the conjecture, it also leads to collapse on the side of identity of proofs.

Prawitz in \cite{P71} found it difficult to justify the left-to-right
direction of the Normalization Conjecture, and Kreisel was looking for
mathematical means that would provide this justification. Maximality is
one such means.

Establishing the left-to-right direction of the Normalization Conjecture
via maximality is like proving the completeness of the classical
propositional calculus with respect to any kind of nontrivial model
via Post completeness (which is proved syntactically by reduction to
conjunctive normal form). Actually, the first proof of this completeness
with respect to tautologies was given by Bernays and Hilbert
exactly in this manner (see \cite{Z99}). 

Maximality for the sort of categories considered above is proved with the
help of coherence in \cite{DP01} and \cite{DP02} (which is established
proof-theoretically, by normalization, cut elimination and similar methods).
Maximality is proved for cartesian closed categories via a typed version
of B\" ohm's theorem in \cite{S83}, \cite{S95} and \cite{DP00}. This
justifies the left-to-right direction of the Normalization Conjecture
also for the implicational and the conjunction-implication fragments
of intuitionistic logic.
The maximality of bicartesian closed categories, which would justify the
left-to-right direction of the Normalization Conjecture for the whole
of intuitionistic propositional logic is, as far as I know, an open problem.
(A use for maximality similar to that propounded here and in \cite{DP00}
and \cite{DP01} was envisaged in \cite{W01}, it seems independently.)

In \cite{D99} (Section 4.11) it is proved that the general notion of
adjunction is also maximal in some sense. The maximality we encountered above,
which involves connectives tied to particular adjunctions, cannot be derived
from the maximality of the general notion of adjunction, but these matters
should not be foreign to each other.

The Normalization Conjecture and the Generality Conjecture do not agree
for the conjunction-implication fragment of intuitionistic logic.
We do not have coherence for cartesian closed categories if the graphs in the
graphical
category are taken to be of the tangle type Kelly and Mac Lane had for
symmetric monoidal closed categories combined with the graphs we had for
cartesian categories (see the preceding section). Both the soundness
part and the completeness part of coherence fail for cartesian closed
categories. The soundness part of coherence also fails for
distributive bicartesian categories, and a fortiori for bicartesian closed
categories. The problem is that in these categories distribution
of conjunction over disjunction is taken to be an isomorphism,
and graphs do not deliver that (see \cite{DP02}, Section 1).

The problem with the soundness part of coherence for cartesian closed
categories may be illustrated with typed lambda terms in the following
manner. By beta conversion and
alpha conversion, we have the following equation:
\[
\lambda_x \langle x,x \rangle \lambda_y y =
\langle \lambda_y y, \lambda_z z \rangle
\]
\noindent for $y$ and $z$ of type $p$, and $x$ of type $p^{\mbox{\it p}}$
(which
corresponds to $p$ {\it implies} $p$).
The closed terms on the two sides of this equation are both of type
$p^{\mbox{\it p}} \times p^{\mbox{\it p}}$.
The type of the term on the left-hand side cannot be
further generalized, but the type of the term
$\langle \lambda_y y, \lambda_z z \rangle$, 
can be generalized to $p^{\mbox{\it p}} \times q^{\mbox{\it q}}$.
The problem noted here does not
depend essentially on the presence of surjective pairing
$\langle \ldots , \ldots \rangle$ and of product types;
it arises also with purely functional types. This problem depends
essentially on the multiple binding of variables, which we have in
$\lambda_x \langle x,x \rangle$; that is, it depends on the structural rule
of contraction.

This throws some doubt on the right-to-left direction of the Normalization
conjecture, which Prawitz, and Kreisel too, found relatively
unproblematic. It might be considered strange that two derivations represent
the same proof if, without changing inference rules, one can be generalized
in a manner in which the other cannot be generalized.

In \cite{F75} (p. 234) Feferman expresses doubt about the right-to-left
direction of the Normalization Conjecture because he finds that a proof
$\pi$ of $A(t)$ ending with the elimination of the universal quantifier
from the premise $\forall xA(x)$ is different from the proof $\pi'$ of
$A(t)$ obtained from $\pi$ by removing $\forall xA(x)$ as a maximum
formula. The proof $\pi$ contains more information than $\pi'$.
In \cite{P81}
Prawitz replied that he takes a proof not ``as a collection of
sentences but as the result of applying certain operations to obtain
a certain end result'' (p. 249). It seems Feferman's objection is
at the level of {\it provability}, and not at the level of
{\it proofs}. A related objection due to Kreisel, and some other
objections to the Normalization Conjecture, may be found in \cite{T75}
(Section 5.3).

The Normalization Conjecture has, for the time being,
the following advantage over the
Generality Conjecture. It applies also to predicate logic, whereas the
latter conjecture has not yet been investigated outside propositional
logic.

When we compare the two conjectures we should also say something about their
computational aspects. With the Normalization Conjecture, 
we have to rely on reduction to a unique normal form in the typed
lambda calculus in order to check equivalence of derivations in the
conjunction-implication fragment of intuitionistic propositional logic.
Nothing more practical than that is known, and such syntactical methods
may be tiresome. Outside of the conjunction-implication fragment, in
the presence of disjunction and negation, such methods become
uncertain.

Methods for checking equivalence of derivations in accordance with the
Generality Conjecture, i.e.\ methods suggested by coherence results, 
often have a clear advantage. This is like the advantage truth tables
have over syntactical methods of reduction to normal form in
order to check tautologicality. However, the semantical methods
delivered by coherence results have this advantage only if the
graphical category is simple enough. And when we enter into categories
suggested by knot theory, this simplicity may be lost. Then, on the
contrary, syntax may help us to decide equality in the graphical
category.

The Normalization Conjecture has made a foray in theoretical
computer science, in the area of typed programming languages. It is not
clear whether one could expect the Generality Conjecture to play a
similar role.

The reflexive and transitive closure of Prawitz's 
immediate reducibility relation, i.e.\ Prawitz's reducibility relation, may
be deemed more important than his equivalence
relation, which we have considered up to now. This matter
leads outside the topic of our survey, which is about
{\it identity} of proofs, but it is worth mentioning. We may ``categorify''
the identity relation between proofs, and consider not only other relations
between proofs, like Prawitz's reducibility relation,
but maps between proofs. The proper framework for doing that
seems to be the framework of weak 2-categories, where we have 2-arrows
between arrows; or we could even go to {\it n}-categories, where we
have {\it n}+1-arrows between {\it n}-arrows (one usually speaks of
{\it cells} in this context). Composition of 1-arrows is associative only
up to
a 2-arrow isomorphism, and analogously for other equations between
1-arrows. Identity of 1-arrows is replaced by 2-arrows satisfying
certain conditions, which are called {\it coherence conditions}.
This notion of coherence is related to the coherence we have
considered above, but it is specific, and need not be the same.
Reducibility between arrows, however, does not
necessarily give rise to an isomorphism.

In the context of the Generality Conjecture, we may also find it natural
to consider 2-arrows instead of identity. The orientation would here be given
by passing from a graph with various ``detours'' to a graph that is
more ``straight'', which need not be taken any more as equal to
the original graph.

With all this we would enter into a very lively field of category theory,
interacting with other disciplines, mainly topology (see \cite{Le01},
\cite{Le02} and papers cited therein). The field looks very promising for
general proof theory, both from Prawitz's and from Lambek's point of
view, but, as far as I know, it has not yet yielded to proof
theory much more than promises.

\section{The Normalization Conjecture in classical logic}

The Normalization Conjecture fares rather well in intuitionistic logic,
but not so well in classical logic. The main difficulty in classical
logic is tied to the fact that in every bicartesian closed category,
for every object $A$ there is at most one arrow with domain $A$ and
the initial object $\bot$ as codomain. Assuming the Normalization
Conjecture, this would mean that in intuitionistic propositional logic
for every formula $A$ there is at most one proof with
assumption $A$ and the absurd constant proposition $\bot$ as conclusion,
i.e.\ there cannot be more than one way to reduce $A$ to the absurd.

In \cite{LS86} the discovery of that fact is credited to Joyal (p. 116),
and the fact is established (on p. 67, Proposition 8.3) by relying on
a proposition of Freyd (see \cite{F72}, p. 7, Proposition 1.12) to the
effect that if in a cartesian closed category {\it Hom}$(A, \bot)$ is
not empty, then $A \cong \bot$; that is, $A$ is isomorphic to $\bot$.
Here is a simpler proof of the same fact.

\vspace{.2cm}

\noindent {\bf Proposition 1.} {\it In every cartesian closed category
with an initial object} $\bot$ {\it we have that} {\it Hom}$(A, \bot)$
{\it is either empty or a singleton.}

\vspace{.2cm}

\noindent {\it Proof.} In every cartesian closed category with $\bot$ we have
$\pi^1_{\bot , \bot}= \pi^2_{\bot , \bot}: \bot \times \bot
\rightarrow \bot$, because {\it Hom}$(\bot \! \times \! \bot, \bot) \cong$
{\it Hom}$(\bot, \bot^{\mbox{$\bot$}})$. Then
for $f,g:A \rightarrow \bot$ we have
$\pi^1_{\bot , \bot} \circ \langle f,g \rangle =
\pi^2_{\bot , \bot} \circ \langle f,g \rangle$, and so $f=g$. {\it q.e.d.}

\vspace{.3cm}

In \cite{LS86} (p. 67) it is concluded from Proposition 1
that if in a bicartesian closed category for every object $A$ we have
$A \cong \lnot \lnot A$,
where the ``negation'' $\lnot B$ is $\bot^{\mbox{\it B}}$
(which corresponds to
$B$ {\it implies} $\bot$), then this category is a preorder. If classical
logic requires $A \cong \lnot \lnot A$ for every proposition $A$,
then the proof theory of that
logic is trivial: there is at most one proof with given assumptions
and a given conclusion.

If the requirement $A \cong \lnot \lnot A$ is deemed too strong, here
is another proposition, which infers triviality from another natural
requirement.

\vspace{.2cm}

\noindent {\bf Proposition 2.} {\it Every cartesian closed category
with an initial object} $\bot$ {\it in which we have a natural
transformation whose components are} $\zeta_A: \lnot \lnot A \rightarrow A$
{\it is a preorder.}

\vspace{.2cm}

\noindent {\it Proof.} Take $f,g: \lnot \lnot A \rightarrow B$, and  
take the canonical arrow 
$\varepsilon_A': A \rightarrow \lnot \lnot A$, which we have by the
cartesian closed structure of our category. Then we have 
$\lnot \lnot (f \circ \varepsilon_A')=\lnot \lnot (g \circ \varepsilon_A')$
by Proposition 1, and from
\[
\zeta_B \circ \lnot \lnot (f \circ \varepsilon_A')=
\zeta_B \circ \lnot \lnot (g \circ \varepsilon_A'),
\]
by the naturality of $\zeta$, we infer
\[
f \circ \varepsilon_A' \circ \zeta_A =  g \circ \varepsilon_A' \circ \zeta_A.
\]
Since $\varepsilon_A' \circ \zeta_A = {\mbox{\bf 1}}_{\lnot \lnot A}$
by Proposition 1, we have $f=g$.

Then, for $\top$ terminal, we have
\[
\begin{array}{rcl}
{\mbox{\it Hom}}(C,D)&\cong &{\mbox{\it Hom}}(\top,D^{\mbox{\it C}})
\\
                     &\cong &{\mbox{\it Hom}}(\lnot \lnot \top,D^{\mbox{\it C}}),
\end{array}
\]
since
$\top \cong \lnot \lnot \top$, and {\it Hom}$(\lnot \lnot \top,D^{\mbox{\it C}})$
is at most a singleton, as we have shown above. {\it q.e.d.}

\vspace{.3cm}

Here is another proposition similar to Proposition 2.

\vspace{.2cm}

\noindent {\bf Proposition 3.} {\it Every bicartesian closed category
in which we have a dinatural
transformation whose components are} $\xi_A: \top \rightarrow A+ \lnot A$
{\it is a preorder.}

\vspace{.2cm}

\noindent {\it Proof.} Take $f,g: \top \rightarrow A$. Then $\lnot f =
\lnot g$ by Proposition 1, and from
\[
({\mbox{\bf 1}}_A + \lnot f) \circ \xi_A =
({\mbox{\bf 1}}_A + \lnot g) \circ \xi_A,
\]
by the dinaturality of $\xi$, we infer
\[
(f+ \lnot {\mbox{\bf 1}}_{\top}) \circ \xi_{\top}=
(g+ \lnot {\mbox{\bf 1}}_{\top}) \circ \xi_{\top}.
\]
Since $\xi_{\top}$ is an isomorphism, we obtain
$f+ \lnot {\mbox{\bf 1}}_{\top}=g+ \lnot {\mbox{\bf 1}}_{\top}$,
from which $f=g$ follows with the help of first injections. {\it q.e.d.}

\vspace{.3cm}

The naturality of $\zeta$ in Proposition 2 is a requirement with
proof-theoretical justification: it has to do with permuting cuts
with rules for negation. Something similar applies to $\xi$.

So, with the Normalization Conjecture, the case for the general
proof theory of
classical logic looks pretty bleak. A number of natural requirements
lead to collapse and triviality.

All these considerations involve an initial object $\bot$. In the
absence of the initiality of $\bot$, matters stand better.
Besides that, classical logic can be embedded by a double-negation
translation
not only in Heyting's intuitionistic logic, but also in
Kolmogorov's minimal intuitionistic logic, which is essentially
Heyting's negationless logic. Relying on this translation, we
could still obtain a nontrivial proof theory for classical logic
in the presence of the Normalization Conjecture.

The initiality of $\bot$, that is the requirement that there is
just one proof with given conclusion and with $\bot$ as assumption,
may perhaps be questioned. If $\bot$ is not devoid of structure,
why could not we proceed in different ways to infer something from it?
For example, if $\bot$ is taken to be $p \wedge \lnot p$, we can infer
$p$ from it either by the first projection, or by passing to
$p \wedge (\lnot p \vee p)$ and then applying the disjunctive
syllogism that enables us to infer $B$ form $A \wedge (\lnot A \vee B)$.
Why must we take these proofs as identical? (We have here two
derivations with different generality in Lambek's sense:
the first generalizes to a derivation of $p$ from $p \wedge \lnot q$ and
the second to a derivation of $q$ from $p \wedge \lnot p$; so,
according to the Generality Conjecture, we have two different proofs.)

Still, rejecting the initiality of $\bot$ looks like a desperate measure,
not in tune with the other intuitions underlying the Normalization Conjecture.
And if we choose to embed classical logic in intuitionistic logic,
classical logic would not be on its own, but would be taken only
as a fragment of intuitionistic logic.

The problems we find with the Normalization Conjecture in
classical logic might accord with the conception
that intuitionistic logic is truly the logic of {\it provability},
whereas classical logic is the logic of something else: it is the logic of
{\it truth} (and falsehood),
and we should not be surprised to find that its general proof theory
is trivial. If, however, we rely on the Generality Conjecture, matters
may look different, as we will see in the next section.

\section{The Generality Conjecture in classical logic}

The Generality Conjecture seems more promising for classical logic than the
Normalization Conjecture. Anyway, it does not lead to collapse and
triviality. The matter has not yet been fully explored, since apart from
fragmentary coherence results, mentioned in Section 4 above, where the
two conjectures agree, and apart from coherence results in substructural
logics, the Generality Conjecture was rather neglected in logic.

The Generality Conjecture underlies the method of proof nets for
multiplicative classical linear logic. The insight provided by
\cite{CS97} is that proof nets are based on a coherence result for
$*$-autonomous categories, which in its turn is based on a coherence
result for categories with the {\it linear distribution}
$A\wedge (B \vee C) \rightarrow (A \wedge B) \vee C$ (previously called
{\it weak distribution}). Linear distribution is what we need for
multiple-conclusion cut. 

The distribution arrow of type
$A\wedge (B \vee C) \rightarrow (A \wedge B) \vee (A \wedge C)$
defined in terms of linear distribution with the help of the
diagonal arrow of type $A \rightarrow A\wedge A$ is not an
isomorphism. This is an important difference with respect to
bicartesian closed categories, where distribution is an isomorphism.

Graphs related to our graphs were tied to sequent derivations of classical
logic in \cite{B91} and \cite{C97}, but apparently without the intention
to discuss the question of identity criteria for proofs.

The price we have to pay for accepting the Generality Conjecture in
classical logic is that not all connectives will be tied to adjoint
functors, as required by Lawvere's Thesis. For example, conjunction and
disjunction will be tied to such functors in the conjunctive-disjunctive
fragment of classical logic, but when we add implication or negation,
in the whole of classical propositional logic, these adjunctions will be
lost. The diagonal arrows of type $A \rightarrow A\wedge A$, for instance,
will not make a natural transformation (this has to do with the
failure of the equation involving $\langle x,x \rangle$, which we
considered in Section 4).

But there might be gains in accepting the Generality Conjecture in
classical logic, and we shall point briefly in the next section
towards a prospect that looks interesting.

\section{Addition of proofs and zero proofs}

Gentzen's multiple-conclusion sequent calculus for classical logic has a
rule of {\it addition} of derivations, which is derived as follows:

$$
\lpravilo{contractions}
\f
{\f
{\f{f:A \rightarrow B}{\theta^R_C(f): A \rightarrow B,C}
{\mbox{\hspace{2em}}}
\f{g:A \rightarrow B}{\theta^L_C(g):C,A \rightarrow B}
}
{{\mbox{\it {cut}}}(\theta^R_C(f),\theta^L_C(g)):A,A \rightarrow B,B}
}
{f+g:A \rightarrow B}
$$

\noindent Here $\theta^R_C(f)$ and $\theta^L_C(g)$ are
obtained from $f$ and $g$
respectively by thinning on the right and thinning on the left,
and {\it cut}$(\theta^R_C(f),\theta^L_C(g))$ may be conceived as 
obtained by applying to $f$ and $g$ a limit case
of Gentzen's multiple-cut rule {\it mix}, where the collection of
mix formulae is empty. (A related rule was
considered under the name {\it mix} in linear logic.)

In a cut-elimination procedure like Gentzen's, $f+g$ is reduced
either to $f$ or to
$g$ (see \cite{G35}, Sections III.3.113.1-2).
If we have $f+g=f$ and $f+g=g$, then we get immediately 
$f=g$, that is collapse and triviality (cf.\ Appendix B.1
by Lafont of \cite{GTL89}). 
If we keep only one of these equations and reject the other,
then to evade collapse we must reject
the commutativity of $+$, but it seems all these decisions would
be arbitrary.
(For similar reasons, even
without assuming the commutativity of $+$, the assumptions of
\cite{S78}, p. 232, C.12, lead to preorder.)  

The Generality Conjecture tells us that we should have neither $f+g=f$ nor
$f+g=g$. The addition of two graphs may well produce a graph
differing from each of the graphs added. It also tells us that addition of
proofs should be associative and commutative.

If we have addition of proofs, it is natural to assume that we also have
for every formula $A$ and every formula $B$ a {\it zero proof}
$\; 0_{A,B}:A \rightarrow B$,
with an empty graph, which with addition of proofs makes the
structure of a commutative monoid. We may envisage having zero proofs
$0_{A,B}:A \rightarrow B$ only for those $A$ and $B$ where there is also
a nonzero proof from $A$ to $B$, but here we consider the more
sweeping assumption involving every $A$ and every $B$.

We should immediately face the complaint that with such zero proofs we have
entered into inconsistency, since everything is provable. That is true,
but not all proofs have been made identical, and we are here not interested
in what is provable, but in what proofs are identical. If it
happens---and with the Generality Conjecture it will happen indeed---that
introducing zero proofs is conservative with respect to
identity of proofs which do not involve zero proofs, then it is
legitimate to introduce zero proofs. Provided it is useful for some
purpose. This is like extending our mathematical theories with
what Hilbert called ``ideal'' objects; like extending the positive integers
with zero, or like extending the reals with imaginary numbers.

What useful purpose might be served by introducing zero proofs? With
addition of proofs, our graphical category in the case of
conjunctive-disjunctive logic turns up to be a category of matrices,
rather than the category of binary relations
(see Sections 3 and 4 above), although the category of binary
relations makes sense too, provided we accept that addition of
proofs is idempotent (i.e.\ $f+f=f$, which means + is {\it union}
rather than addition of proofs). Composition becomes matrix
multiplication, and addition is matrix addition. And in the presence of
zero matrices, we obtain a unique normal form like in linear algebra:
every matrix is the sum of matrices with a single 1 entry.

A number of logicians have sought a link between logic and linear algebra,
and here is such a link. We have it not for an alternative logic, but
for classical logic. We have it, however, not at the level of
provability, but at the level of identity of proofs.

The unique normal form suggested by linear algebra is not unrelated to
cut elimination. In the graphical category of matrices, cut elimination is
just matrix multiplication. And the equations of this category yield a
cut-elimination procedure. They yield it even in the absence of zero proofs
(provided $f+f=f$),
and unlike cut-elimination procedures for classical logic
in the style of Gentzen, the new procedure admits a commutative addition of
proofs without collapse. So, in classical logic, the Generality Conjecture
is not foreign to cut elimination, and it would not be foreign to
the Normalization Conjecture if we understand the equivalence
relation involved in this conjecture in a manner different from Prawitz's.

This need not exhaust the advantages of having zero proofs. They may be
used also to analyze disjunction elimination.  Without pursuing this
topic very far, let us note that passing from $A\vee B$ to $A$ involves
a zero proof from $B$ to $A$, and passing from $A\vee B$ to $B$ involves
a zero proof from $A$ to $B$. If next we are able to reach $C$ both from
$A$ and from $B$, we may add our two proofs from $A\vee B$ to $C$, and
so to speak ``cancel'' the two zero proofs.

On an intuitive level, we should not imagine that zero proofs are
faulty proofs. They are rather hypothetical proofs, postponed proofs,
or something like the oracles of recursion theory.

\section{Propositional identity}

Once we have answered the question of identity criteria for proofs, by
the Normalization Conjecture, the Generality Conjecture,
or in another satisfactory manner, armed with our answer we may
try approaching other
logical questions, as envisaged by Prawitz in his topic 2.4 of
general proof theory (see the beginning of Section 2 above).
Here we will sketch how we can apply this answer 
to settle the question of
propositional identity (which was mentioned by Prawitz in topic 2.3).

Categorial proof theory
suggests concepts that were not previously
envisaged in logic. Such a concept is
isomorphism between sentences, which is understood
as in category theory. (A certain school of category
theory recommends that we should drop identity and
think always in terms of isomorphism.)
A sentence $A$ is {\it isomorphic} to a sentence $B$
if and only if there is a proof $f$ from $A$ to $B$ and a
proof $f^{-1}$ from $B$ to
$A$ such that $f$ composed with $f^{-1}$
is equal to the identity proof
from $A$ to $A$ and $f^{-1}$ composed with $f$ is equal 
to the identity proof from $B$ to $B$.
An identity proof is such that when composed with any
other proof, either on the side of the assumptions or on the side
of the conclusion, it is equal to this other proof.
That two sentences are isomorphic means that they behave
exactly in the same manner in proofs: by composing, we can always extend
proofs involving one of them, either as assumption or as conclusion, to
proofs involving the other, so that nothing is lost, nor gained. There
is always a way back. By composing further with the inverses, we return to
the original proofs.

Isomorphism between sentences is an equivalence
relation stronger than mutual implication. So, for example, $A \wedge B$
is isomorphic to $B \wedge A$,
while $A \wedge A$ only implies and is
implied by $A$, but is not isomorphic to it. (The problem is that the
composition of the proof from $A \wedge A$ to $A$, which is
either the first or the second projection, with
the diagonal proof from $A$ to $A \wedge A$
is not equal to the identity proof from $A \wedge A$ to $A \wedge A$.)

It seems reasonable to suppose that isomorphism between sentences
analyzes propositional
identity, i.e.\ identity of meaning for sentences:

\begin{quotation}
\noindent Two sentences express the
same proposition if and only if they are isomorphic.
\end{quotation}

\noindent This way we would base propositional identity
upon identity of proofs, since in defining isomorphism of
sentences we relied essentially on a notion of identity of proofs.

By relying on either the Normalization Conjecture or the
Generality Conjecture, the formulae isomorphic in conjunctive
logic are characterized by the equations of commutative monoids.
(By duality, this solves the problem also for disjunctive logic.)
By relying on the Normalization Conjecture,
the formulae isomorphic in the conjunction-implication fragment of
intuitionistic logic, i.e.\ objects isomorphic in all cartesian closed
categories, have been characterized in \cite{S81} via the
axiomatization of equations with multiplication, exponentiation and one,
true for natural numbers. This fragment of arithmetic (upon which one comes
in connection with Tarski's ``high-school algebra problem'') is finitely
axiomatizable and decidable.
(For an analogous result about the multiplicative conjunction-implication
fragment of intuitionistic linear logic see \cite{DP97}.)
The problem of characterizing isomorphic formulae 
in the whole intuitionistic propositional calculus, which
corresponds to bicartesian closed categories, seems to be still open.

Leibniz's analysis of identity, given by the equivalence

\begin{quotation}
\noindent Two individual terms name the same object if and only if
in every sentence one can be replaced by the other without change
of truth value,
\end{quotation}

\noindent assumes as given and
uncontroversial propositional equivalence, i.e.\ identity of truth value for
sentences, and
analyzes in terms of it identity of individuals.
We analyzed above  
propositional identity in terms of isomorphism between sentences,
a notion that presupposes an understanding of identity of proofs, and
our analysis resembles Leibniz's to a certain extent.

Propositional equivalence, which in classical logic is defined by identity
of truth value, amounts to mutual implication, and
is understood as follows in a proof-theoretical context:

\begin{quotation}
\noindent $A$ is equivalent to $B$ if and only if
there is a proof from $A$ to $B$
and a proof from $B$ to $A$.
\end{quotation}

\noindent This relation between the propositions $A$ and $B$,
which certainly cannot
amount to the stricter relation of propositional identity, does not 
presuppose understanding identity of proofs.
We analyzed propositional identity 
by imposing on the composition of the two proofs
from $A$ to $B$ and from $B$ to $A$ a condition formulated in terms
of identity of proofs.

The question of propositional identity seems as philosophical as the
question of identity of proofs, if not more so, but the answer we proposed
to this question is rather mathematical, as were the answers given
to the question of identity of proofs by the Normalization
and Generality Conjectures.

\section{Conclusion}

The question we have discussed here suggests a perspective in logic---or
perhaps
we may say a dimension---that has not been explored enough.
Logicians were, and still are, interested mostly in provability, and
not in proofs. This is so even in proof theory. When we address the
question of identity of proofs we have certainly left the realm of
provability, and entered into the realm of proofs. This should become clear
in particular when we introduce the zero proofs 
of Section 7.

The complaint might be voiced that with the Normalization and Generality
Conjectures we are giving very limited answers to the question of identity of
proofs. What about identity of proofs in the rest of mathematics, outside
logic? Shouldn't we take into account many other inference rules, and not
only those based on logical constants? Perhaps not if the structure of proofs
is taken to be purely logical. Perhaps conjectures like the Normalization
and Generality Conjectures are not far from the end of the road.

Faced with two concrete proofs in mathematics---for example, two
proofs of the Theorem of Pythagoras, or something more
involved---it could seem pretty hopeless
to try to decide whether they are identical
just armed with the Normalization Conjecture or the Generality Conjecture.
But this hopelessness might just be the hopelessness of formalization.
We are overwhelmed not by complicated principles, but by sheer quantity.

An answer to the question of identity criteria for proofs could perhaps
also be expected to come from the area of proof theory tied to
complexity theory, which we mentioned at the end of the first section.
There is perhaps an equivalence relation between derivations defined
somehow in terms of complexity by means of which we could
formulate a conjecture analogous to the Normalization and Generality
Conjectures, but I am not aware that a conjecture of that kind
has been proposed.

Although the problem of identity criteria for proofs is a conceptual
mathematical problem, one might get interested in it for
technical reasons too. There is interesting mathematics
behind this problem. But one should probably also have an interest in logic
specifically. Now, it happens that many logicians are not interested
in logic. They are more concerned by other things in mathematics.
They dwell somehow within the realm of what is accepted as logic,
but they deal with other mathematical matters. Logic is taken for granted.
Some logicians, however, hope for discoveries in logic.

\vspace{.5cm}

\noindent {\footnotesize {\it Acknowledgements.} I would like to thank
warmly Zoran Petri\' c for a very great number of conversations, in
which ideas exposed here were reached or pursued.
A plenary address at the Logic Colloquium
2002 in M\" unster was based on the present paper.
I am grateful to the organizers for inviting me
to deliver this talk. I am also grateful to
the Alexander von Humboldt Foundation for supporting my participation
at the colloquium. The writing of the paper was financed by
the Ministry of Science, Technology and Development of
Serbia through grant 1630 (Representation of proofs with applications,
classification of structures and infinite combinatorics).}

\end{document}